\documentclass{amsart}
\usepackage{color}

\usepackage{color}

\usepackage{amssymb}
\usepackage{amsmath}
\usepackage{pdfsync}

\def\CX{{\mathbb C}}
\def\QX{{\mathbb Q}}
\def\NX{{\mathbb N}}

\def\ZX{{\mathbb Z}}
\def\PX{{\mathbb P}}

\def\GL{{\rm GL}}

\def\gl{{\rm gl}}

\def\diag{{\rm Diag}}
\def\d{{
\partial}}

\def\calS{{\mathcal S}}

\def\calF{{\mathcal  F}}

\def\calK{{\mathcal K}}

\def\calN{{\mathcal  N}}

\def\calD{{\mathcal D}}
\def\calC{{\mathcal C}}
\def\calF{{\mathcal F}}
\def\calO{{\mathcal O}}
\def\calI{{\mathcal I}}
\def\calOU{{\calO_{\calU}}}
\def\calOUU{{\calO_{\calU'}}}

\def\calS{{\mathcal S}}
\def\calU{{\mathcal U}}
\def\calV{{\mathcal V}}

\def\d{{
\partial}}
\def\dx{{
\partial_{x}}}

\def\Fhat{{\hat{\calF}}}
\def\Nhat{{\hat{\nabla}}}
\def\pp{{\mathbb{P}}^1({\mathbb{C}})}

\def\diag{{\rm diag}}
\def\Ga{{{\mathbb G}_a}}
\def\Gm{{{\mathbb G}_m}}

\theoremstyle{definition}
\newtheorem{theorem}{Theorem}[section]
\newtheorem{lemma}[theorem]{Lemma}

\newtheorem{definition}[theorem]{Definition}
\newtheorem{example}[theorem]{Example}

\newtheorem{lem}[theorem]{Lemma}
\newtheorem{cor}[theorem]{Corollary}
\newtheorem{prop}[theorem]{Proposition}

\numberwithin{equation}{section}

\begin{document}

\title[Monodromy groups  of 
parameterized  differential equations]{Monodromy groups  of \\ 
parameterized linear differential equations with regular singularities}


\author{Claude Mitschi}
\address{Institut de Recherche Math\'ematique Avanc\'ee, Universit\'e de Strasbourg et CNRS, 7 rue Ren\'e Descartes, 67084 Strasbourg Cedex, France}
\curraddr{}
\email{mitschi@math.unistra.fr}

\author{Michael F. Singer}
\address{Department of Mathematics, North Carolina
State University, Box 8205, Raleigh, North Carolina 27695-8205}
\curraddr{}
\email{singer@math.ncsu.edu}
\thanks{The second author was partially supported by NSF Grants CCF-0634123 and CCF-1017217. }


\date{}

\dedicatory{}

\commby{}

\begin{abstract}

We study the notion of regular singularities for parameterized complex ordinary linear differential systems, prove  an analogue of the Schlesinger theorem for systems with  regular singularities and  solve both  a parameterized version of the weak Riemann-Hilbert Problem and a special case of the inverse problem in parameterized Picard-Vessiot theory.
\end{abstract}
\maketitle
\section {Introduction}{  Let 
\begin{eqnarray}\label{INTRO:eqn1}
\frac{dY}{dx} &= &A(x) Y
\end{eqnarray}
be a linear differential equation where $A(x) \in \gl_n(\CX(x))$ is an $n\times n$ matrix with coefficients that are rational functions over the complex numbers. One can associate two groups to such an equation,  the {monodromy group} and the {differential Galois group}.  To define the monodromy group one starts by removing
the set $S =\{a_1, \ldots, a_s\}$ of singular points (possibly including infinity) of (\ref{INTRO:eqn1})  from the Riemann sphere $\pp$ and fixing a point $a_0 \in \pp\backslash S$. Using standard existence theorems,  there exists a fundamental solution matrix, that is, an $n\times n$ matrix $Z = (z_{i,j})$ of functions analytic in a neighborhood of $a_0$ with $\det Z(a_0) \neq 0$. Analytic continuation of $Z$ along any closed  path $\gamma$ in $\pp\backslash S$ centered at $a_0$ yields a new fundamental
solution matrix $Z_\gamma$ which is related to $Z$ via an equation of the form $Z_\gamma =ZM_\gamma$ for some $M_\gamma \in \GL_n(\CX)$. One can show that $M_\gamma$  depends only on the homotopy { class } of $\gamma$ in $\PX^1(\CX)\backslash S$ and that the map $\gamma \mapsto M_\gamma$ defines a homomorphism  $\rho:\pi_1(\PX^1(\CX)\backslash S; a_0) \rightarrow \GL_n(\CX)$, that is, a representation of  the fundamental group. }{ Selecting a different fundamental solution matrix results in a conjugation of the image of $\rho$. The image of $\rho$ is called the {\em monodromy group} of (\ref{INTRO:eqn1}) and is determined up to conjugacy.   \\[0.1in]
To define the differential Galois group of (\ref{INTRO:eqn1}), one forms the field $K =$\linebreak $ \CX(x, z_{1,1}, \ldots , z_{n,n})$ constructed from $\CX(x)$ by adjoining the entries of $Z$. Note that (\ref{INTRO:eqn1}) implies that this field is closed under the action of the derivation $\frac{d}{dx}$. The { differential Galois group }, also called the {\em Picard-Vessiot group} $G$ is the group of all field-theoretic automorphisms of $K$ which leave any element of $\CX(x)$ fixed and commute with~$\frac{d}{dx}$ (see \cite{DAAG} or \cite{PuSi2003} for an exposition of the associated theory). One can show that for any $\sigma \in G$, $\sigma(Z) = (\sigma(z_{i,j})) = Z M_\sigma$ for some $M_\sigma \in \GL_n(\CX)$. From the preservation of algebraic relations under analytic continuation, one can show that the {\em monodromy group} $\rho(\pi_1(\PX^1(\CX)\backslash S; a_0))$ is contained in the Picard-Vessiot group~$G$.  The Picard-Vessiot group  is closed in the Zariski topology, that is, there exists a system of polynomial equations in $n^2$ variables such that~$G$ is precisely the set of invertible matrices whose entries satisfy these equations. Furthermore, various properties of solutions of (\ref{INTRO:eqn1}) are reflected in properties of $G$.  For example the dimension of $G$ (as an algebraic or complex Lie group) is related to the algebraic dependence among the $z_{i,j}$,  and the solvability of the connected component containing the identity of $G$ is equivalent to (\ref{INTRO:eqn1}) being solvable in terms of exponentials, integrals and algebraic functions. \\[0.1in]
When one restricts the type of singular points of (\ref{INTRO:eqn1}), one can say more concerning the relationship of the monodromy group and the { Picard-Vessiot group}.  We say a singular point $a \in S$ is  {\em regular singular}  if there is an $n\times n$ matrix $P(x)$ of functions { meromorphic} at $a$ with $\det P(a) \neq 0$ such that 
the matrix $U = PZ$ satisfies an equation of the form
\begin{eqnarray}\label{INTRO:eqn2}
\frac{dU}{dx} & = & \frac{\tilde{A}}{x-a} U
\end{eqnarray}
where $\tilde{A}$ is a {\em constant} matrix, i.e., $\tilde{A} \in \gl_n(\CX)$ (there is an equivalent definition of regular singular in terms of the growth of the entries of $Z$ near $a$. See Chapters 3.1 and 5.1 of \cite{PuSi2003} for a fuller discussion). A result of Schlesinger ({\cite{schlesinger}, \S\,159, 160;} \cite{PuSi2003}, Theorem 5.8) states that if  all the singular points of (\ref{INTRO:eqn1}) are regular singular, then the  Picard-Vessiot group is the smallest Zariski-closed subgroup of $\GL_n(\CX)$ containing the monodromy group. \\[0.1in]
%
One may also consider inverse questions, that is, which groups appear as monodromy or { Picard-Vessiot groups}. For example, one may ask: given a homomorphism $\rho:\pi_1(\PX^1(\CX)\backslash S; a_0) \rightarrow \GL_n(\CX)$, does there exist an equation (\ref{INTRO:eqn1}) whose monodromy group is the image of $\rho$? This is a version of the so-called {\em Riemann-Hilbert Problem} and has a positive solution (see Chapters 5 and 6 of \cite{PuSi2003} for a fuller discussion of the various versions of this problem as well as other references).  Using a solution of this problem, C.~and M.~Tretkoff showed that any Zariski closed subgroup of $\GL_n(\CX)$ is the { Picard-Vessiot} group of some equation (\ref{INTRO:eqn1}) over $\CX(x)$.\\[0.1in]
In this paper we consider similar results for parameterized systems of linear differential equations. Parameterized families of linear differential systems with regular singular points arise in the study of isomonodromic as well as monodromy evolving deformations and their relation to the equations of mathematical physics (\cite{Bol_iso_def}, \cite{ChAb,ChAb2}, \cite{its}, \cite{IKSY}, \cite{sabbah}, \cite{ohyama,ohyama2}, \cite{MS}) .  We address some analogous fundamental questions concerning  the monodromy groups of such { families. { More precisely, we consider parameterized linear differential systems of the form
\begin{eqnarray} \label{leteqpar}\frac{\d Y}{\d x} &=& A(x,t)Y\end{eqnarray}}
{where the entries of the matrix $A(x,t)$ are rational functions of $x$ with coefficients that are analytic in the multiparameter $t$ on some domain of $\CX^r$.} We begin, in Section 2, by studying equivalent definitions of  regular singular points of  such systems and  proving bounds on the growth of solutions in the neigborhood of these singularities. }{  In Section 3, we  
  show that 
the parameterized monodromy matrices of    
a system \eqref{leteqpar} belong to its associated parameterized Picard-Vessiot group. In Section 4, we prove an analogue of the Schlesinger theorem for systems with  regular singularities. Our result states that  for such systems    the parameterized monodromy matrices generate a Kolchin-dense subgroup  of the parameterized Picard-Vessiot group. In Section 5, we  solve a parameterized version of the Riemann-Hilbert Problem   and a special case of the inverse problem in parameterized Picard-Vessiot theory.}}

\section{Parameterized regular singularities} 
Let $\calU$ be an open connected  subset of $\CX^r$ with  $0\in \calU$, and let $\calO_{\calU}$ be the ring of analytic functions on $\calU$ of a variable $t$. Let $\alpha\in \calOU$ and assume $\alpha(0)=0$. We will denote  
\begin{enumerate}
\item by $\calOU((x-\alpha(t)))$ the ring of formal Laurent series in powers of $x-\alpha(t)$ with coefficients in $\calOU$, that is, elements  
\[  f(x,t)=\sum_{i\ge m} a_i(t)(x-\alpha(t))^i\]
where $m\in \ZX$ is independent of $t$,

\item by $\calOU(\{x-\alpha(t) \})$  the ring  of  those $f(x,t)\in\calOU((x-\alpha(t)))$  that,   for  each fixed  $t\in \calU$, converge  for $0<|x-\alpha(t)|<R_t$, for some $R_t>0$.
\end{enumerate} 

\begin{lemma}\label{radius}
Let $f(x,t)\in\calOU(\{x-\alpha(t) \})$ and let    $\calN\subset\calU$ be a compact neighborhood of $0$.  Then there is $R>0$ such that the series  $f(x,t)$ converges for all $t\in\calN$ and $0<|x-\alpha(t)|<R$.
 \end{lemma}

 \begin{proof} For each $t\in \calU$  we may assume that $R_t$  is maximal, possibly infinite. For finite $R_t$,  let  $\Gamma(\alpha(t),R_t)$ and $D(\alpha(t),R_t)$  denote the circle and open disk respectively,  with center $\alpha(t)$ and radius $R_t$.   If $R_{t_0}=\infty$ for some $t_0\in \calU$ then clearly $R_t=\infty$ for all $t\in \calU$. Assuming this is not the case, $R_t$ is a continuous function of~$t$. To prove this, fix $t_0\in \calU$, and  a neighborhood $u(t_0)$ of $t_0$ in $\calU$ such that $\alpha(t)\in D(\alpha(t_0),R_{t_0})$ for all $t\in u(t_0)$. If $t\in u(t_0)$,  the circles  
$\Gamma(\alpha(t_0),R_{t_0})$ and $\Gamma(\alpha(t),R_t)$ either are equal, or intersect at two points, or are inner tangent. A simple geometric argument shows that
\[ |R(t)-R(t_0)|\le |\alpha(t)-\alpha(t_0) | \]
for all $t\in u(t_0)$, and the continuity of $R$ follows from the continuity of $\alpha$. Since $\calN$ is compact, and $R_t>0$ for all $t\in \calN$, the function $R_t$ (possibly infinite) has a lower bound $R>0$ on $\calN$.
\end{proof}

Consider a parameterized linear differential equation
\begin{eqnarray}\label{system}
{\frac{\d Y}{\d x}}&=&A Y
\end{eqnarray}
where $A\in\gl_n\bigl(\calOU(\{x-\alpha(t)\})\bigr)$. Note that $A$  may be writen as
\[
A(x,t)=\frac{A_{-m}(t)}{(x-\alpha(t))^{m}}+\frac{A_{-m+1}(t)}{(x-\alpha(t))^{m-1}}+\ldots=\sum_{i\ge{-m}} {(x-\alpha(t))^{i}}{A_i(t)}
\]
where $A_i(t)\in\gl_n\bigl(\calOU\bigr)$ for all $i\ge{-m}$, and $m\in\NX$ does not depend on~$t$.

\begin{definition}\label{equiv}

Two equations
\[ {\frac{\d Y}{\d x}}=AY    \quad{\mathrm{and}}\quad    {\frac{\d Y}{\d x}}=BY, \]
with $A,B\in\gl_n\bigl(\calOU(\{x-\alpha(t)\})\bigr)$, are {\em equivalent} if there exists $P\in \GL_n\bigl(\calOU(\{x-\alpha(t)\})\bigr)$ such that
\[B={\frac{\d P}{\d x}}P^{-1}+PAP^{-1}, \]
that is, if $Y$ satisfies the first equation, then $PY$ satisfies the second. 
\end{definition}

\begin{definition}\label{regsing} With notation as before,

\begin{enumerate}
\item Equation (\ref{system}) has {\em simple singular points  near $0$} if $m=1$ and $A_{-1}\ne 0$ as an element of $\gl_n\bigl(\calOU(\{x-\alpha(t)\})\bigr)$,

\item Equation (\ref{system}) has {\em parameterized regular singular points near $0$} if it is equivalent to an equation with simple singular points near $0$.
\end{enumerate}
\end{definition}
{ Note that in the non-parameterized case, simple singular points are sometimes referred to as ``Fuchsian singular points'' and regular singular points are sometimes referred to as ``regular points''. 

\begin{example} Let 
\begin{eqnarray*}
A & = &\left(\begin{array}{cc} 0 & -3\\ 0& 0 \end{array}\right) \frac{1}{ (x-t)^2} + \left(\begin{array}{cc} t & 0\\ 0& t-2 \end{array}\right) \frac{1}{ x-t}\\
B & = & \left(\begin{array}{cc} t-1&0\\ 0& t-1 \end{array}\right) \frac{1}{x-t}
\end{eqnarray*}
A calculation shows that $B = \frac{\d P}{\d x} P^{-1} + PAP^{-1}$ where 
\begin{eqnarray*}
P & = & \left(\begin{array}{cc} \frac{1}{x-t}&\frac{-1}{(x-t)^2}\\ 0& x-t \end{array}\right) 
\end{eqnarray*}
Therefore,  { the equations $\frac{\d Y}{\d x} = AY$ and $\frac{\d Y}{\d x} = BY$ are equivalent. Since $\frac{\d Y}{\d x} = BY$ has  simple singular points near $0$, the equation $\frac{\d Y}{\d x} = AY$ has parameterized regular singular points}   near $0$.\\[0.1in]
\end{example}

In the previous example, we transformed an equation with  regular singular points near $0$ into an equation that not only has simple singular points but is of the form }{ $\frac{\d Y}{\d x} = \frac{\tilde{A}(t)}{x-\alpha(t)}Y$.  }{ We shall now show that this can be done in general.  }
Let $\delta=(x-\alpha(t))\frac{\d}{ \d x}$. If  Equation (\ref{system}) has simple singular points near $0$ we may also write it as
\begin{eqnarray}\label{fuchs}
 \delta Y&=&\Bigl(\sum_{i\ge 0} (x-\alpha(t))^i A_i(t)\Bigr)Y
\end{eqnarray}
with  (renamed) $A_i\in\gl_n\bigl(\calOU(\{x-\alpha(t)\})\bigr)$  and $A_0\ne 0$. We will show that such an equation is equivalent to an equation of a simpler form. The proof is a slight modification of the similar one for non-parameterized equations in \cite{PuSi2003}.

\begin{prop}\label{normtransf} 
Assume   that in Equation  (\ref{fuchs})
no eigenvalues of $A_0(0)$ differ by positive integers. Then there is an open connected subset $\calU'$  of $\calU$ and matrices $P_i\in\gl_n(\calOUU)$ such that 
\begin{itemize}
\item[1)] the  substitution $Z=PY$, with  
\[ P(x,t)=I + \sum_{i\ge 1}\left(x-\alpha(t)\right)^i P_i(t),\]
transforms Equation (\ref{fuchs}) into
\begin{eqnarray*}
\delta Z &=& A_0(t)Z,
\end{eqnarray*}
\item[2)] the series $P(x,t)$ converges for $(x,t)\in D(\alpha(t),R)\times \calU'$, for some $R>0$ that does not depend on $t$. 
\end{itemize}
\end{prop}
\begin{proof}
We begin by noting that since no eigenvalues of $A_0(0)$ differ by positive integers, there is a neighborhood $\calU'\subset\calU$ of $0$ such that for $t\in\calU'$, no eigenvalues of $A_0(t)$ differ by positive integers. We will now  follow the proof of Proposition 3.12 of (\cite{PuSi2003}), p. 64). We wish to construct a matrix $P$, 
\[ P(x,t)=I+\sum_{i\ge 1}\left( x-\alpha(t)\right)^i P_i(t)\]
such that \[A_0(t)P(x,t)=P(x,t)A(x,t)+\delta P(x,t).\]
Comparing powers of $x-\alpha(t)$ we see that
\begin{equation}\label{gaugeqn} A_0 P_i-P_i\left(A_0+i I \right)= A_i + A_{i-1}P_1 + \ldots + A_1 P_{i-1}
\end{equation}
for all $i\ge 1$ (with $P_0=I$). As noted above, for fixed $t\in\calU'$ no eigenvalues of $A_0(t)$ differ by positive integers. Therefore,  for $t\in\calU'$, the map
\[X\mapsto A_0(t)X-X\left( A_0(t)+iI\right) \]
is an  isomorphism on $\gl_n(\CX)$ and the matrix $M(t)$ representing this $\CX$-linear map has a nonzero determinant. This implies that  the entries of  $M(t)^{-1}$  are analytic on $\calU'$ and that we can solve  Equation (\ref{gaugeqn}) to find matrices $P_i(t)$ whose entries are analytic on $\calU'$. 

We now turn to the statement concerning convergence. The formal  power series $P(x,t)$ (in powers of $x-\alpha(t)$) satisfies the differential equation 
\[\delta P(x,t)=A_0(t)P(x,t)-P(x,t)A(x,t). \]
For each fixed value of $t$ this is a differential equation with a simple (Fuchsian) singularity. Lemma 3.9.2 of \cite{Sibuya} or the proof of Lemma 3.42 of \cite{PuSi2003} implies that $P(x,t)$ has a radius of convergence at least as large as that for $A(x,t)$. We may assume   that $\overline{\calU'}$ is compact with $\overline{\calU'}\subset\calU$ and hence, 	by Lemma \ref{radius} above, that there is  an $R>0$ such that $A(x,t)$ converges for all $t\in\calU'$ and $0<|x-\alpha(t)|<R$, which ends the proof of 2).
\end{proof}
Let us deduce from Proposition \ref{normtransf} a slightly weaker result without the hypothesis   on $A_0(0)$.

\begin{cor}\label{normfuchs} Consider  the general equation (\ref{fuchs}). Then there exists
\begin{itemize}
\item  a constant matrix $C\in \GL_n(\CX)$,

\item an $n\times n$ matrix 
\begin{eqnarray*}
S=\left( \begin{array}{cccc}
(x-\alpha(t))^{r_1} I_1 & 0 & \ldots & 0 \\
0   & (x-\alpha(t))^{r_2} I_2 & \ldots & 0 \\
\vdots   & \vdots  & \vdots & \vdots \\
0   & 0 & \ldots & (x-\alpha(t))^{r_s} I_s  \\
\end{array}
\right)
\end{eqnarray*}
where for some $s$  the $r_i$, $i=1,\ldots, s$, are nonnegative integers and the $I_i$ are identity matrices of various sizes, 

\item an open connected subset $\calU'$  of $\calU$  and matrices $P_i\in\gl_n(\calOUU)$, 
\end{itemize}
such that
 \begin{itemize}
 
\item[1)] the substitution $Z=PSCY$, with  
\[ P(x,t)=I + \sum_{i\ge 1}(x-\alpha(t))^i P_i(t),\]
transforms (\ref{fuchs}) into
\begin{eqnarray}\label{zeqn} \delta Z&=&{\tilde{A}}(t) Z
\end{eqnarray} 
where ${\tilde{A}} \in\gl_n(\calOUU)$,

\item[2)] the series $P(x,t)$ converges for $(x,t)\in D(\alpha(t),R)\times \calU'$ for some $R>0$ that does not depend on $t$. 

 \end{itemize}

\end{cor}
\begin{proof}
We will perform a {\it shearing transformation}  to replace Equation (\ref{fuchs}) by an equation satisfying the condition of Proposition \ref{normtransf} on $A_0(0)$, as follows.  Let $c_1,\ldots,c_s$ be the distinct eigenvalues of  $A_0(0)$ and assume that $c_2-c_1=m$, a positive integer. Replacing, if needed, $Y$ by $CY$ for some $C\in\GL_n(\CX)$, we may assume that
\begin{eqnarray*}
A_0(0)=\left( \begin{array}{cccc}
A_{0,1} & 0 & \ldots & 0 \\
0   & A_{0,2}  & \ldots & 0 \\
\vdots   & \vdots  & \vdots & \vdots \\
0   & 0 & \ldots & A_{0,s}   \\
\end{array}
\right)
\end{eqnarray*}
where $A_{0,i}$, for $i=1,\ldots,s,$ is a matrix in Jordan normal form with eigenvalue $c_i$. Let
\begin{eqnarray*}
T=\left( \begin{array}{cccc}
(x-\alpha(t))^m I_1 & 0 & \ldots & 0 \\
0   &  I_2 & \ldots & 0 \\
\vdots   & \vdots  & \vdots & \vdots \\
0   & 0 & \ldots &I_s  \\
\end{array}
\right)
\end{eqnarray*}
where each $I_i$ is an identity matrix of the same size as $A_{0,i}$. The substitution  $Z=TY$ transforms Equation (\ref{fuchs}) into an equation 
\begin{eqnarray*}
\delta Z=\Big(\sum_{i\ge 0} (x-\alpha(t))^i {\tilde{A}}_i\Big)Z
\end{eqnarray*}
where the eigenvalues of ${\tilde{A}}_0(0)$  are $c_1+m,c_2,\ldots,c_s$, that is, ${\tilde{A}}_0(0)$ has fewer eigenvalues differing by positive integers. By induction, one constructs a matrix $S$ such that the transform of Equation (\ref{fuchs}) via  $Z=SCY$ satisfies the condition of Proposition \ref{normtransf}, whose conclusion ends the proof.  
\end{proof}
The following corollary shows that our definition of  a parameterized regular singularity yields, as in the non-parameterized case, solutions that have moderate growth in the neighborhood of this moving singularity.

\begin{cor}\label{modgrowth} Assume that Equation (\ref{system}) has regular singular points near $0$. Then there is an open  connected subset $\calU'$  of $\calU$ such that
\begin{itemize}
\item[1)]  Equation (\ref{system}) has a solution $Y$ of the form
\begin{eqnarray}\label{sol}
Y(x,t)&=&\Bigl( \sum_{i\ge i_0}(x-\alpha(t))^i Q_i(t) \Bigr) (x-\alpha(t))^{{\tilde{A}}(t)}
\end{eqnarray}
with  ${\tilde{A}}\in\gl_n(\calOUU)$ and $Q_i\in\gl_n(\calOUU)$ for all $i\ge i_0$,

\item[2)] for any $r$-tuple $(m_1,\ldots,m_r)$ of nonnegative integers there is an integer $N$ such that for any fixed $t\in{\mathcal U}'$ and  any sector ${\mathcal S}_t$ from $\alpha(t)$ in the complex plane, of opening less than $2\pi$, 
\[Ê\lim_{x\rightarrow\alpha(t)\atop x\in {\mathcal S}_t} \bigl(x-\alpha(t)\bigr)^N \frac{\d ^{m_1+\ldots+m_r}Y(x,t)}{\d^{m_1}t_1\ldots\d^{m_r}t_r}=0.\] 

  \end{itemize}

\end{cor}

\begin{proof}
Under our asumptions, Equation (\ref{system})  is equivalent to an equation with  simple singularities near $0$. We shall prove the conclusion for this new equation and one easily sees that it  holds for the original equation as well. 

Assume Equation (\ref{system}) has  simple singular points near $0$. Corollary \ref{normfuchs} implies that there are $C$, $S$ and $P$ as described { above} such that  $Y$ is a solution of  Equation (\ref{fuchs})  whenever   $Z=PSCY$ is a solution of Equation (\ref{zeqn}). This shows that
\[ Y=(PSC)^{-1} \bigl(x-\alpha(t) \bigr)^{{\tilde{A}}(t)}\]
is a solution of the desired form (\ref{sol}). Differentiating (\ref{sol}) yields a form that satisfies the conclusion of 2). 
\end{proof}

Solutions of parameterized differential equations with irregular singularities have been studied in  \cite{babbitt_varadarajan} and \cite{schaefke}. Assuming $0$ is a (non-moving)  irregular singularity, these authors gave a  condition on the exponential part of a formal solution in the usual form
\[ {\hat{Y}}(z)={\hat{H}}(z)z^Je^Q\]
to ensure that the coefficients of the formal series ${\hat{H}}(z)$ depend analytically on the multiparameter.
\section{Parameterized monodromy}

\subsection{Classical Picard-Vessiot theory and monodromy}
Consider a differential equation
\begin{eqnarray*} \frac{dY}{dx}& = &A(x)Y\end{eqnarray*}
where $A \in \gl_n(\QX(x))$. To apply differential Galois theory to this equation we need to work over an algebraically closed field containing $\QX$, for example $\bar{\QX}$, the algebraic closure of $\QX$. But when we talk about monodromy matrices and want to say that the monodromy matrices lie in the  Picard-Vessiot group,  we need to account for the possibility that these matrices have  transcendental entries. To deal with this, we usually  consider the Picard-Vessiot theory over $\CX(x)$, but we might use any algebraically closed field   containing the entries of the monodromy matrices.\\[0.1in]
This raises the question: does the group change when we go to a bigger field of constants?  The answer is given by the following proposition.  We note that Propositions~\ref{letprop1} and~\ref{letprop2} below can be stated in far greater generality than stated here (see \cite{GGO}, Theorem 9.10 for an approach via Tannakian categories and also \cite{CHS07}, pp.~80-81 for similar results concerning difference equations) yet we present a proof here using simple tools that allows us to prove this result in our restricted setting. \\[0.1in]
 If not otherwise specified, the derivation will be denoted by $(\ )'$. We will write PV-group and PV-extension, for short,  for the Picard-Vessiot group and  Picard-Vessiot extension respectively. 

\begin{prop}\label{letprop1} 
Let $C_0\subset C_1$ be algebraically closed fields and $k_0 = C_0(x), k_1 = C_1(x)$ be differential fields where $c' = 0$ for all $c \in C_1$ and $x'=1$.
Let 
\begin{eqnarray}\label{leteqn0}Y' &=& AY\end{eqnarray}
be a differential equation with $A
 \in \gl_n(k_0)$. If $G(C_0)\subset\GL_n(C_0)$ is the PV-group over $k_0$ of Equation (\ref{leteqn0}) with respect to some fundamental solution,  where $G$ is a linear algebraic group defined over $C_0,$  then $G(C_1)$ is the PV-group of (\ref{leteqn0}) over $k_1$,  with respect to some fundamental solution. 
\end{prop}
\begin{proof}
Equation (\ref{leteqn0}) has a regular point in $C_0$ and we shall assume this is $x=0$. We can think of $C_i(x)$ as a subfield of the field (of formal Laurent series) $C_i((x))$ for $i=0,1$. Since $0$ is a regular point, Equation (\ref{leteqn0}) has a solution $Z_0 \in \GL_n(C_0[[ x]])\subset\GL_n(C_0((x)))$ (note that this solution is found by substituting $Z_0 = Z_{0,0} + Z_{0,1}x +Z_{0,2}x^2 + \ldots$ into (\ref{leteqn0}) and equating powers of $x$. Assuming $Z_{0,0} = I_n$ ensures that $Z_0 \in \GL_n(C_0[[x]])$. Since the differential field  $K_0 = k_0(Z_0) \subset k_0((x))$ has no new constants, it  is a PV-extension  of $k_0$ for Equation (\ref{leteqn0}). Similarly, $K_1 = k_1(Z_0)$ is  a PV-extension  of $k_1$ for Equation (\ref{leteqn0}).

Let $R_0 = k_0[Z_0]$ and $R_1 = k_1[Z_0]$.  We may write $R_0 = k_0[Y]/\calI_0$ and $R_1= k_1[Y]/\calI_1$ where $Y$ is a set of $n^2$ indeterminates and $\calI_0$ and $\calI_1$ are ideals in their respective rings. Let us prove that $\calI_1 = k_1\cdot \calI_0$.  Clearly, $k_1\cdot \calI_0 \subset \calI_1$. To show the other inclusion, let $P \in \calI_1$.  By clearing denominators, we may assume that $P \in C_1[x][Y]$.  Let us write $P = \sum d_i P_i$ where $\{d_i\}$ is a $C_0$-basis of $C_1$ and $P_i \in C_0[x][Y]$. Note that the $d_i$ are linearly independent over $C_0[[x]]$ as well. Substituting $Z_0$ for $Y$ in $P$, we  have that $P(Z_0) = \sum d_i P_i(Z_0)=0$. Therefore $P_i \in \calI_0$ for each $i$ and so $P \in k_1\cdot \calI_0$.

  For $i = 0,1$, the PV-group of $K_i$ over $k_i$ consists of the matrices $B$ in $\GL_n(C_i)$ such that  $Z_1B$ is again a zero of $\calI_i$.   Let $\calI_0 = (f_1, \ldots , f_m)$ be generated in $k_0[Y]$ by polynomials $ \ f_j \in k_0[Y] $.  The PV-group of $K_0$ over $k_0$ then consists  of all $B \in \GL_n(C_0)$ such that $f_j(Z_0B) = 0$ for $j = 1, \ldots ,m$.  We may consider $f_j(Z_0B)$ as an element of $C_0((x))$ and write $f_j(Z_0B) = \sum_{\ell \ge \ell_0} f_{\ell,j}(B) x^\ell$ where $f_{\ell,j}(B)$ is a polynomial in the entries of $B$  with coefficients  {\it in $C_0$}. The PV-group of $K_0$ over $k_0$ is therefore defined as 
\[G_0=\{B \in \GL_n(C_0) \ |\  f_{\ell,j}(B)=0, \mbox{ for  }j=1,\ldots,m \mbox{ and }  \ell\ge \ell_0, \}.\]
 Since $\calI_1 = k_1\cdot \calI_0$, we  have $\calI_1=(f_1, \ldots , f_m)$  in $k_1[Y]$, and the    PV-group of $K_1$ over $k_1$ is also defined as  $G_1=\{ B \in \GL_n(C_1) \ | \ f_j(Z_0B) = 0 \mbox{ for } j = 1, \ldots ,m\}$ and therefore as
\[G_1=\{B \in \GL_n(C_1) |  f_{\ell,j}(B)=0, \mbox{ for  } j=1,\ldots,m \mbox{ and }  \ell\ge \ell_0. \}.\] 
  This means that if $G$ is a linear algebraic group defined over $C_0$ by $\{f_{\ell,j}(B)=0\}$  then the PV-group of $K_i$ over $k_i$ is $G(C_i)$ for $i=0,1$. To end the proof, note that the linear algebraic group $G$ defined over $C_0$ is uniquely determined by the group $G(C_0)$ of its $C_0$-points.
  \end{proof}

\begin{cor}
 Assume  in Equation (\ref{leteqn0}) that $A \in \gl_n(C_0(x))$ where $C_0$ is some algebraically closed subfield of $\CX$. Assuming  $0$ is a regular point, let us  fix it as the base-point of $\pi_1(\PX^1(\CX)\backslash \calS)$, where  $\calS$ is the set of singular points of (\ref{leteqn0}) on  $\PX^1(\CX)$. Let $G(C_0)$ be the PV-group of (\ref{leteqn0}) over $C_0(x)$, where $G$ is a linear algebraic group defined over $C_0$. If $C_1$  is any algebraically closed subfield of $\CX$  containing $C_0$ and the entries of the monodromy matrices, then the monodromy matrices are elements of the PV-group $G(C_1)$ of  (\ref{leteqn0}) over $C_1(x)$.
\end{cor}

\begin{proof}
With notation from the proof of Proposition~\ref{letprop1}, note that the matrix $Z_0$ has entries that are {\it convergent} in some neighborhood of $0$ since at a regular point formal solutions are convergent. Let $M$ be a monodromy matrix with respect to $Z_0$  corresponding to an element $[\gamma]$ of $\pi_1(\PX^1(\CX)\backslash \calS)$ and let $P \in \calI_1$, that is, $P(Z_0) = 0$.  We again may assume that $P \in C_1[x][Y]$. To show that $M \in G(C_1)$ it is enough to show that $P(Z_0M) = 0$.  Note that we can analytically continue $P(Z_0)$ around $\gamma$ and the result will be $P(Z_0M)$. Since $P(Z_0) = 0$, we have $P(Z_0M) = 0$ as well.
\end{proof}

\medskip
\subsection{Parameterized Picard-Vessiot theory and monodromy}
In the pa\-ra\-me\-teri\-zed Picard-Vessiot theory (PPV-theory for short)  one has a similar issue: the equation may have coefficients that lie in one differentially closed field while the parameterized monodromy matrices have entries that lie in a  larger differentially closed field. For general definitions and facts about PPV-theory we refer to \cite{CaSi} and \cite{Landesman}.

We will first prove a result similar to Proposition~\ref{letprop1} for parameterized Picard-Vessiot extensions (PPV-extensions for short). 
For simplicity of notation, we shall consider equations of the form
\begin{eqnarray} \label{leteqn1}\frac{\d Y}{\d x} &=& A(x,t)Y\end{eqnarray}
where $A(x,t)\in\gl_n(\calOU(x))$, { with $t=(t_1,\ldots,t_r)$} in  some domain $\calU\in\CX^r$. We shall denote differentiation with respect to $x, t_1,\ldots,t_r$ by $\d_x, \d_{t_1},\ldots,\d_{t_r}$ respectively, and let $\Delta=\{\d_x, \d_{t_1},\ldots,\d_{t_r}\}$  and $\Delta_t
=\{\d_{t_1},\ldots,\d_{t_r}\}$. 
Let $C$ be a differentially closed  $\Delta_t$-extension
 of some field of functions that are  analytic  on some domain of~$\CX^r$ and let $\d_{t_i}$ denote for each $i$  the derivation extending $\d_{t_i}$. We define the differential   $\Delta$-field $k =C(x)$,  where $x$ is an indeterminate over $C$,  by  letting $\d_x(x) = 1, \d_{t_i}(x) = 0$ for each $i$ and $\d_x(c) = 0$ for all $c\in C$. We willl always assume that $C$ is chosen such that $A\in gl_n(k)$. For the definition of a differentially closed field, see (\cite{CaSi}, Definition 3.2).

PPV-extensions  of $k$ for Equation (\ref{leteqn1}) are of the form $K=k\langle Z\rangle$ where $K$ has no new $\d_x$-constants and $Z$ is a fundamental solution for (\ref{leteqn1}). The brackets $\langle\ldots \rangle$ denote the fact that $K$ is generated,   as a $\Delta_t$-field,  by the entries of $Z$, that is, $K = k( Z, \d_{t_1}Z, \ldots , \d_{t_r}Z, \ldots, \d_{t_1}^{\alpha_1}\ldots \d_{t_r}^{\alpha_r}Z,\ldots)$. Note that this is automatically a $\d_x$-field as well, with $\d_xZ = AZ, \ \d_x(\d_{t_i}Z) = (\d_{t_i}A)Z + A \d_{t_i}Z$ for each $i$,  etc.
\begin{prop}\label{letprop2}Let $C_0\subset C_1$ be differentially closed $\Delta_t$-fields and let $k_0 = C_0(x)$,  $k_1 = C_1(x)$ be $\Delta$-fields as above. Let 
\begin{eqnarray}\label{leteqn2}\d_x Y&=& AY\end{eqnarray}
be a differential equation with $A \in \gl_n(k_0)$.  If $G(C_0)\subset\GL_n(C_0)$ is the PPV-group over $k_0$ of Equation (\ref{leteqn2}) with respect to some fundamental solution,  where $G$ is a linear differential algebraic group defined over the differential $\Delta_t$-field $C_0,$  then $G(C_1)$ is the PPV-group over $k_1$ of (\ref{leteqn2}) with respect to some fundamental solution. 
\end{prop}

\begin{proof} The proof is basically the same as in the non-parameterized case. Equation (\ref{leteqn2}) has a regular point in $C_0$ and we shall assume this is $x=0$. We can think of $C_i(x)$ as a subfield of the field (of formal Laurent series) $C_i((x))$ for $i=0,1$. Since $0$ is a regular point, Equation (\ref{leteqn2}) has a solution $Z_0 \in \gl_n(C_0((x)))$.  Assuming (using the previous notation $Z_0=\sum Z_{0,i}x^i$) that $Z_{0,0} = I_n$ ensures that $Z_0 \in \GL_n(C_0((x)))$. Since $K_0 = k_0\langle Z_0\rangle \subset k_0((x))$ has no new $\d_x$-constants, it is a PPV-extension  of $k_0$ for Equation (\ref{leteqn2}). Similarly $K_1 = k_1\langle Z_0\rangle$ is  a PPV-extension  of $k_1$ for (\ref{leteqn2}). \

Let $R_0 = k_0\{Z_0\}$ and $R_1 = k_1\{Z_0\}$ (where $\{\ldots \}$ denotes the {\it differential ring} generated by $\ldots$).  We may write $R_0 = k_0\{Y\}/\calI_0$ and $R_1= k_1\{Y\}/\calI_1$ where $Y$ is a set of $n^2$ {\it differential} indeterminates and $\calI_0$ and $\calI_1$ are  {\it differential} ideals in their respective rings. Exactly as  in the non-parameterized case, one shows that $\calI_1 = k_1\cdot \calI_0$. 


 For $i = 0,1$,  the PPV-group of $K_i$ over $k_i$ is  the set of matrices $B$ in $\GL_n(C_i)$ such that  $Z_0B$ is again a zero of $\calI_i$.   Let $\calI_0 = \{f_j\}_{j \in J}$ for some indexing set $J$.  We then have that  the PPV-group of $K_0$ over $k_0$ is the set of $B \in \GL_n(C_0)$ such that $f_j(Z_0B) = 0$ for $j \in J$.  We may consider $f_j(Z_0B)$ as an element of $C_0((x))$ and write $f_j(Z_0B) = \sum_{\ell \geq \ell_0} f_{\ell,j}(B) x^\ell$ where $f_{\ell,j}(B)$ is a {\it differential} polynomial in the entries of $B$  with coefficients {\it in $C_0$}. Therefore the PPV-group of $K_0$ over $k_0$ is defined as
\[G_0=\{B \in \GL_n(C_0) | f_{\ell,j}(B)=0, \mbox{ for  }j\in J \mbox{ and }  \ell\ge \ell_0, \}.\]
 Since $\calI_1 = k_1\cdot \calI_0$, the differential ideal $\calI_1$ is generated by $\{f_j\}_{j \in J}$ in $k_1\{Y\}$, and hence the PPV-group of $K_1$ over $k_1$ is defined 
 as 
 \[G_1=\{B \in \GL_n(C_1) | f_{\ell,j}(B)=0, \mbox{ for  }j\in J \mbox{ and }  \ell\ge \ell_0, \}.\]
This means that if $G$ is the linear differential  algebraic group defined over $C_0$ by $\{f_{\ell,j}(B)=0\}$,   then the PPV-group of $K_i$ over $k_i$ is $G(C_i)$ for $i=0,1$, and since $G$ is defined over $C_0$ by the group of its $C_0$-points, this ends the proof. 
\end{proof}

We shall now show that the parameterized monodromy matrices (defined below) of Equation (\ref{leteqn1}) are elements of the PPV-group. \\[0.1in]
Let $\calD$ be an open subset of $\pp$ with $0 \in \calD$.   Assume that $\pp\backslash \calD$ is the union of $m$ disjoint disks $D_i$ and that for each $t \in \calU$, Equation (\ref{leteqn1})  has a unique singular point in each $D_i$.  Let $\gamma_i$,  $i = 1, \ldots , m$ be the obvious  loops  generating $\pi_1(\calD, 0)$. Let us fix a fundamental solution~$Z_0$ of (\ref{leteqn1}) in the neighborhood of $0$ and define, for each fixed $t\in\calU$, the monodromy matrices of (\ref{leteqn1}) with respect to this solution and the $\gamma_i$. We will call these matrices, which depend on $t$, the {\it parameterized monodromy matrices} of Equation  (\ref{leteqn1}). \\[0.1in]
%
%
 In the classical, non-parameterized situation, we gave an argument using analytic continuation of $P(Z_0)$ where $P$ had coefficients in $C_0(x)$. This argument made sense because the coefficients of $P$ were {\it a fortiori} analytic functions. In the present situation the coefficients, {\it a priori}, do not have such a meaning. The following result of Seidenberg \cite{sei58,sei69} allows us to give them such a meaning.
%

\begin{theorem}[Seidenberg]\label{seidenberg} Let $\QX\subset \calK \subset \calK_1$  be finitely generated differential extensions of the field of rational numbers $\QX$, 
and assume that  $\calK$   consists of
meromorphic functions on some domain $\Omega\in\CX^r$. Then  $\calK_1$ is isomorphic to a field $\calK_1^*$ of functions that are meromorphic on a domain $\Omega_1 \subset \Omega$, such that $\calK|_{\Omega_1} \subset \calK_1^*$.
\end{theorem}
%
%
\noindent Note that finitely generated here means finitely generated in the differential sense. We are now able to prove the same result as in the PV-case.

\begin{theorem}\label{mongal}
Assume in Equation (\ref{leteqn1}) that $A\in\gl_n(C_0(x))$, where   $C_0$ is any differentially closed $\Delta_t$-field containing  $\CX$ and let  $C_1$ be  any differentially closed $\Delta_t$-field containing $C_0$ and the entries of the parameterized monodromy matrices  of Equation (\ref{leteqn1}) with respect to a fundamental solution of (\ref{leteqn1}). Then the parameterized  monodromy matrices belong to   $G(C_1)$, where $G$ is the PPV-group of (\ref{leteqn1}) over the  $\Delta$-field $C_0(x)$. 
\end{theorem}

\begin{proof}
As in the proof of Proposition \ref{letprop2}, let us construct a formal fundamental solution $Z_0 \in \GL_n(C_0((x)))$ of Equation (\ref{leteqn2}). The coefficients of the power series $Z_0$ are matrices whose entries are polynomials, with integer coefficients, in the $x$-coefficients of the entries of $A(x,t)$, assuming as we are that $Z_0$ has the initial condition $Z_0(t,0)=I_n$. The coefficients of $Z_0$ are therefore analytic on some domain~$\calU\subset \CX^r$. The usual estimates show that $Z_0$ is analytic in $\calV\times \calU$ where $\calV$ is a neighborhood of $0$ in the $x$-plane. Let $\calK$ be the differential $\Delta_t$-field generated over $\QX$ by the $x$-coefficients of the entries of $A(x,t)$ (that is, the coefficients of powers of $x$ of these  rational functions). Note that $\calK$ consists of  functions meromorphic on $\calU$. Let $M$ be a parameterized monodromy matrix corresponding to some $[\gamma]$. For any differential polynomial $P \in C_0(x)\{Y\}$ let us show that $P(Z_0M)=0$ whenever $P(Z_0)=0$. As in the PV-case, this will prove that $M\in G(C_0)$.  Let $\calK_1$ be a finitely generated $\Delta_t$-field extension of $\calK$ such that $M \in \GL_n(\calK_1)$ and $P \in \calK_1(x)\{Y\}$.. The Seidenberg theorem applied to $\calK$ and $\calK_1$ allows us to identify $\calK_1$ with a differential field of meromorphic functions on some $\calU_1\subset\calU$ and in particular to consider $P(Z_0)$ and $P(Z_0M)$  as  functions analytic on $\calV\times \calU_2$ for some $\calU_2\subset\calU_1$. Let us continue $P(Z_0)$ along $\gamma$, considered for each fixed $t\in\calU_2$ as a path in $\calD\times\{t\}\subset\calD\times\calU_2$. Since the coefficients of $P$ remain unchanged,  $P(Z_0M)=0$ holds after continuation along $\gamma$.  \end{proof}
\section{A parameterized analogue of the Schlesinger theorem}\label{schlsec}

We consider as before a parameterized equation
\begin{eqnarray}\label{leteqn3}\frac{\d Y}{\d x} &=& A(x,t)Y\end{eqnarray}
where $A(x,t)$ is a rational function of $x$ with coefficients that are functions of a multi-variable $t=(t_1,\ldots,t_r)$, analytic in some domain $\calU\in\CX^r$ containing $0$. We shall denote differentiation with respect to $x$ and $t_1,\ldots,t_r$ by $\d_x$ and $\d_{t_1},\ldots,\d_{t_r}$ respectively and let as before $\Delta =\{\d_x,\d_{t_1},\ldots,\d_{t_r}\}$,  $\Delta_t=\{\d_{t_1},\ldots,\d_{t_r}\}$. 

The following lemma is inspired by a result of R. Palais (\cite{Pa78}). 

\begin{lem}\label{palais} Let $\calF$ be a $\Delta$-field of functions that are meromorphic on $\calV \times \calU$ where $\calV\subset\CX$ and $\calU\subset\CX^r$ are open connected sets, and let $\calC_x=\{u\in\calF\  |\  \d_xu=0 \}$. Furthermore assume $x \in \calF$. Let $f\in\calF$ be such that $f(x,t)\in\CX(x)$ for each $t\in\calU$. Then for some $m\in\NX$, there exist $a_0,\ldots,a_m,b_0,\ldots, b_m\in\calC_x$ such that
\[ f(x,t)=\frac{\sum_{i=0}^m a_ix^i}{\sum_{i=0}^m b_ix^i}\]
\end{lem}

\begin{proof}
For each $r\in\NX$, let
\[X_r=\{t\in\calU\ |\ f(x,t)\ {\mathrm{is\ a\ ratio\ of\  polynomials\ of\ degrees\ at\ most\ } r} \}. \]
By asumption $\calU=\cup_{r\in\NX} X_r$. The Baire Category Theorem implies that for some $m$, the closure of  $X_m$ has a nonempty interior. For $t\in X_m$, there exist 
$a_{0,t}, \ldots,a_{m,t}$, $b_{0,t},\ldots, b_{m,t} \in\CX$ such that
\[Ê f(x,t)=\frac{\sum_{i=0}^m a_{i,t}x^i}{\sum_{i=0}^m b_{i,t}x^i}.\]Ê
This implies that, for $t\in X_m$,
\[ x^mf(x,t), x^{m-1}f(x,t), \ldots,f(x,t),\ x^m,\ldots,1\] 
are linearly dependent over $\CX$. In particular the Wronskian determinant (with respect to $\d_x$) 
$W(x,t)=wr(x^mf(x,t), x^{m-1}f(x,t), \ldots,f(x,t),x^m,\ldots,1)$ vanishes for any $t\in X_m$. Since the closure of $X_m$ has a nonempty interior, $W(x,t)=0$ on a nonempty open set and so is identically zero. Since \[x^mf(x,t), x^{m-1}f(x,t), \ldots,f(x,t),\ x^m,\ldots,1 \in\calF\] 
the vanishing of the Wronskian determinant implies that these elements are also linearly dependent over the $\d_x$-constants of $\calF$, that is, $\calC_x$. Therefore there exist  $a_0,\ldots,a_m,b_0,\ldots, b_m\in\calC_x$, not all zero, such that
$(\sum_{i=0}^m b_ix^i)f(x,t)=\sum_{i=0}^m a_ix^i$. 
\end{proof}
\noindent This lemma is the key  to proving an analogue of the Schlesinger theorem. 
 
 \smallskip
\noindent As in the previous section, let $\calD$ be an open subset of $\pp$ with $0 \in \calD$, assuming that $\pp\backslash \calD$ is the union of $m$ disjoint disks $D_i$ and that for each $t \in \calU$, Equation (\ref{leteqn3})  has a unique singular point $x_i(t)$ in each $D_i$.  Let $\gamma_i, i = 1, \ldots , m$ be the obvious  loops  generating $\pi_1(\calD, 0)$. We also assume that $0\in\calU$. We construct as before a fundamental solution $Z_0$ of (\ref{leteqn3}) near $0\in\calD$, analytic on $\calV\times \calU$ where $\calV$ is some neighborhood  of $0$ in $\calD$, and such that $Z_0(t,0)=I_n$. By the theorem about analytic dependence on initial conditions ({\it cf.} \cite{Cartan}), the analytic continuation of $Z_0$ along each $\gamma_i$ (for each fixed $t\in\calU$) provides a solution which is again analytic on $\calV\times\calU$. We may therefore assume, on a possibly smaller $\calU$,  that $B_{\gamma_i}\in\GL_n(\calOU)$ for each parameterized monodromy matrix $B_{\gamma_i}$ with respect to $Z_0$. 

\smallskip
In what follows $k(x)$, for any differential $\Delta_t$-field~$k$, will denote the field of rational functions in the indeterminate~$x$ and coefficients in~$k$, where $x$ is  a $\Delta_t$-constant and $\dx$ is  the usual, formal differentiation of rational functions with $\dx(x)=1$ and~$\dx(a)=0$ for all~$a\in k$.

 \begin{theorem}\label{schlesinger}
With notation as before, assume that Equation (\ref{leteqn3}) has regular singularities near each $x_i(0)$, $i=1,\ldots,m$. Let $k$ be a  differentially closed $\Delta_t$-field containing  the $x$-coefficients of the entries of $A$, the singularities $x_i(t)$ of (\ref{leteqn3}) and the entries of the parameterized monodromy matrices with respect to $Z_0$. Then the parameterized monodromy matrices  generate a Kolchin-dense subgroup of $G(k)$, where $G$ is the PPV-group of (\ref{leteqn3})  over $k(x)$.
\end{theorem}

\begin{proof}
Under the asumptions of the theorem,   $K=k(x)\langle Z_0\rangle$ is a PPV-extension of $K$  for (\ref{leteqn3}) (it has no new constants) and by Theorem \ref{mongal} we know that the parameterized monodromy matrices lie in its PPV-group $G(k)$. We now wish to show that these generate a Kolchin-dense subgroup.  Using the parameterized Galois correspondence it is enough to show,  for any $f\in K$, that $f\in k(x)$ whenever $f$ is left invariant by the action of the monodromy matrices. Let $\calF_0$ be the differential $\Delta_t$-subfield of $k$ generated over $\QX$ by the coefficients of powers of $x$ in the entries of~$A$, the singularities~$x_i(t)$  and the entries of the parameterized monodromy matrices with respect to~$Z_0$. Note that $\calF_0$ consists of  functions  meromorphic on $\calU$. Let us fix $f\in K$, and any  finitely generated differential $\Delta_t$-subfield $\calF_1$ of $k$ containing  $\calF_0$,  such that $f\in\calF_1(x)\langle Z_0\rangle$. Using Theorem \ref{seidenberg} we may consider $\calF_1$ as a field of meromorphic functions on some (possibly smaller than $\calV\times\calU$) domain $\Omega$ of the $(x,t)$-space. We may moreover assume that $f$ is a rational function, with coefficients that are analytic functions on $\Omega$ (rational in $x$, analytic in $t$) in the entries  of both $Z_0$ and  some $\Delta_t$-derivatives of $Z_0$. Assume that $f$ is left invariant by the parameterized monodromy matrices, that is, by analytic continuation of $f$ along each $\gamma_i$, for any fixed value of $t$ such that $(0,t)\in\Omega$. For such fixed values of  $t$, Corollary \ref{modgrowth} implies that  $f$ has moderate growth at its singular points. Therefore, for these fixed values of $t$, $f$ is a rational function of $x$ with coefficients in $\CX$. Lemma \ref{palais} implies that $f$ is a rational function of $x$ with coefficients in the subfield $\calC_x$ of $\d_x$-constants of $\calF_1(x)\langle Z_0\rangle$, and since $\calC_x\subset k$, this ends the proof.
\end{proof}

\section{A weak  parameterized Riemann-Hilbert Problem} \noindent Classically, the  weak form of the Riemann-Hilbert Problem is:

\smallskip
\begin{quotation}\noindent {\em  Let $S=\{a_1, \ldots a_s\}$ be a finite subseet  of $\PX^1(\CX)$ and $a_0\in\PX^1(\CX)\backslash S$. Given a representation  \[\rho: \pi_1(\PX^1(\CX)\backslash S; a_0) \rightarrow \GL_n(\CX)\]
of the fundamental group $\pi_1(\PX^1(\CX)\backslash S; a_0)$ show that there exists a linear differential system 
\[\frac{dY}{dx} = A Y\] with $A\in \gl_n(\CX(x))$, having only regular singular points, all  in $S$,  such that for some  fundamental solution analytic at $a_0$, the  monodromy representation is   $\rho$.}\end{quotation}

\smallskip

 \noindent Solutions of this problem go back to Plemelj with modern versions presented by R\"ohrl, Deligne and others (see \cite{anosov_bolibruch} or \cite{PuSi2003} for presentations of solutions and historical references). In this section, we will present a solution to a parameterized version of this problem and apply this to the inverse problem in parameterized Picard-Vessiot Theory. \\[0.1in]
\noindent To state a parameterized version of the weak Riemann-Hilbert Problem, note that a representation of $\pi_1(\PX^1(\CX)\backslash S ; a_0)$ above is determined by the images of the generators of this group, that is, by selecting $s$ invertible matrices $M_1, \ldots , M_s$ such that $M_1\cdot\ldots\cdot M_s = I_n$, the identity matrix. We will present a solution of the following weak parameterized Riemann-Hilbert Problem

\smallskip

 \begin{quotation}\noindent {\em  Let $S=\{a_1, \ldots a_s\}$ be a finite subset  of $\PX^1(\CX)$ and $a_0\in \PX^1(\CX)\backslash S$, and let  $D$ be an open  polydisk\footnote{An open polydisk $D$ is a set of the form $\{t = (t_1, \ldots , t_r)   \ | \ |t_j - w_j | < r_j \ 1\leq j \leq r\}$ for some $w = (w_1, \ldots , w_r)$ and positive real numbers $r_j$.} in $\CX^r$. Let $\gamma_1, \ldots , \gamma_s$ be generators of $\pi_1(\PX^1(\CX)\backslash S ; a_0)$ and, for $i =1 , \ldots , s$, let $M_i:D\rightarrow \GL_n(\CX)$ be analytic maps with $M_1\cdot\ldots\cdot M_s=I_n$. Show that there exists a parameterized linear differential system 
 \begin{eqnarray*}\label{eqnRH1}
 \frac{\d Y}{\d x} &= &A(x,t)Y
 \end{eqnarray*}
with $A\in \gl_n(\calO_{D'}(x))$ for some open polydisk $D' \subset D$,  having only regular singular points, all in $S$,  and such that for some fundamental solution, the parameterized monodromy matrix along each $\gamma_i$ is $M_i$.}
\end{quotation}

\smallskip
\noindent Let $K$ be the field of functions, meromorphic on $D$, generated over $\QX$ by the entries of the matrices $M_i$.  The assignment $\gamma_i \mapsto M_i$ yields a homomorphism $\chi:\pi_1(\PX^1(\CX)\backslash S ; a_0)\rightarrow \GL_n(K)$. \\[0.1in]
Modern solutions of the weak Riemann-Hilbert Problem use the techniques of analytic vector bundles and we will proceed in a similar fashion following the presentation in \cite{BoMaMi}. We begin by first constructing a  family of vector bundles on $\PX^1(\CX)$ together with a family of meromorphic connections.  We will then show how this yields a solution of the weak parameterized Riemann-Hilbert Problem.{ {We note that the related but different problem of constructing an {\it isomonodromic} family including a given Fuchsian differential equation is considered in \cite{miw81} and \cite{ma80a} (see also \cite{Heu} and \cite{sabbah}).}} \\[0.1in]
Vector bundles are determined by cocycles with respect to coverings so we will proceed by defining parameterized cocycles on $\PX^1(\CX)$. 
 We begin  by considering $\PX^1(\CX)\backslash S$ (and will fill in the ``holes'' later). \\[0.1in]
Let $U_{s+1}, \ldots , U_N$ be a covering of 
$\PX^1(\CX)\backslash S$ by open disks (other disks $U_1, \ldots ,U_s$ will be defined below) and let $\eta_i$ be a path in $\PX^1(\CX)\backslash S$ from $x_0$ to the center of $U_i$. For each pair  $(i,j)$ with $i<j$ and $U_i\cap U_j \neq \emptyset$, let $\delta_{i,j}$ be a the line-segment joining the center of $U_i$ to the center of $U_j$. Let
\[g_{i,j} = \chi([\eta_i\circ \delta_{i,j}\circ\eta_j^{-1}])  \in \GL_n(K).\]
 By replacing $D$ with a smaller open polydisk  if necessary, we may assume that  the entries of all the  $g_{i,j}$ are analytic functions on $D$. For each $t \in D$, we consider $g_{i,j}(t): U_i\cap U_j \rightarrow \GL_n(\CX)$ as a constant function on $U_i\cap U_j$.  One can easily check that, for fixed $t\in D$, the $\{g_{i,j}(t)\}$ form a cocycle and therefore define an analytic vector bundle $\Fhat(t)$ of rank $n$ over $\PX^1(\CX)\backslash S$. For later use it is important to note that the the $\{g_{i,j}\}$ can also be thought of as forming a cocycle  on $(\PX^1(\CX)\backslash S) \times D$ with respect to the covering $\{U_i\times D\}$ and so define an analytic vector bundle $\Fhat$ on $(\PX^1(\CX)\backslash S) \times D$.\\[0.1in]
For each $t \in D$, we define a system of linear differential equations locally over each $U_i, i=s+1, \ldots, N$ via the forms
\[dy = \omega_i(t) y\]
where each $\omega_i(t) = 0$. Of course, this is just the trivial system and one has, on each nonempty intersection $U_i\cap U_j$,
\[ \omega_i(t) = d(g_{i,j}(t))(g_{i,j}(t))^{-1} + (g_{i,j}(t))\omega_j(t)(g_{i,j}(t))^{-1}\]
since the $g_{i,j}(t)$ are constant  on these sets. Therefore these local systems patch together to form an analytic connection $\Nhat(t)$ on $\Fhat(t)$. Trivially, if we write $\Nhat(t)$ in terms of local coordinates, one has that the terms appearing depend analytically on $t$ (and this holds for any other trivializing covering as well). Furthermore, for each $t \in D$, the monodromy associated with simple loops $\gamma_i$ from $x_0$ surrounding each $x_i$ corresponds to $M_i(t)$.\\[0.1in]
We will now fill in the ``holes'' and extend  each  $(\Fhat(t),\Nhat(t))$ to a holomorphic vector bundle $\calF(t)$ on all of $\PX^1(\CX)$ and meromorphic connection $\nabla(t)$ on $\calF(t)$. Let $U_1, \ldots , U_s$ be  pairwise disjoint open disks centered at $a_1, \ldots , a_s$ (for simplicity, we will assume that these are all finite points).
Shrinking $D$ again if necessary (but keeping the same notation), there exists  for each $i$ an analytic  function $N_i(t):D \rightarrow \gl_n(\CX)$ such that $N_i(t) = (1/2\pi i)\log (M_i(t))$. On each $U_i$, $1\leq i\leq s$,  consider for fixed $t \in D$ the meromorphic  system 
\begin{eqnarray}\label{logsys}
dy &= &\frac{N_i(t)}{(x-a_i)}d(x-a_i) \ y.
\end{eqnarray}
For each $i$, $1\leq i \leq s$, select an $\alpha$ such that $U_i\cap U_\alpha \neq \emptyset $ and let $Y_i(x,t)$ be a solution of (\ref{logsys}) on this latter set. We may write 
\[Y_i(x,t) = (x-a_i)^{N_i(t)}.\]
on $U_i\cap U_\alpha$.  Let $g_{i,\alpha} = Y_i(x,t)$. If $U_\beta, \beta\neq \alpha,$  also has nonempty intersection with $U_i$, we consider a path in $U_i$ from a designated point $u_i \in U_i\cap U_\alpha$ ending in $U_i\cap U_\beta$ moving in the counterclockwise direction (less than one turn) around $a_i$. We let $g_{i,\beta}$ denote the analytic continuation of $g_{i,\alpha}$ along this path. One can show that 
 the $g_{k,l}$ thus defined  for all $1\le k,l\le N$ such that $U_k\cap U_l\ne\emptyset$, define a cocycle   for the covering $\{ U_l\times D\}_{1\leq l\leq N}$ of $\PX^1(\CX) \times D$.  This cocycle yields an analytic vector bundle $\calF$ on $\PX^{1}(\CX)\times D$.  {\em A fortiori}, for each $t \in D$, $\{g_{k,l}(t)\}$ defines a cocycle for the covering $\{ U_l\}_{1\leq l\leq N}$ of $\PX^1(\CX)$ and this yields, for each $t \in D$ an analytic vector bundle $\calF(t)$ on $\PX^1(\CX)$. \\[0.1in]
We now claim that for each $t\in D$ the local systems $\{dy = \omega_i(t) y\}_{i=1}^{N}$ define a meromorphic connection on $\PX^{1}(\CX)$. We need only check compatibility on sets of the form $U_i\cap U_\alpha$ with $1\leq i \leq s, s+1 \leq \alpha \leq N$.  On such a set we have
\[d(g_{i,\alpha}(t))g_{i,\alpha}(t)^{-1} + (g_{i,\alpha}(t))\omega_\alpha(t) (g_{i,\alpha}(t))^{-1}= dY_iY_i^{-1} = \omega_i(t)\]
which proves compatibility.  Therefore, for each $t \in D$,  we have a connection $\nabla(t)$ on $\calF(t)$. Note that in local coordinate, the terms of $\nabla(t)$ depend analytically on~$t$.  We  refer to the pair $(\calF(t), \nabla(t))$ as the  {\em canonical extension} of $(\Fhat(t),\Nhat(t))$, and this generalizes similar notions in the non-parameterized case introduced by Deligne. Note that the connection $\nabla(t)$ has at worst logarithmic poles at the  $a_i, i=1, \ldots , s$. Therefore,  the differential systems will have regular (even Fuchsian) singular points in the  local coordinates.  \\[0.1in]
We now wish to change the covering of $\PX^1(\CX)\times D$ and consider the vector bundle $\calF$ with respect to this new covering.  In particular, let $p_1 \neq p_2$ be points of $\PX^1(\CX)$ and let  $V_1 = \PX^1(\CX)\backslash \{p_2\}$ and $V_2 = \PX^1(\CX)\backslash \{p_1\}$. We claim that $\calF$ is isomorphic to a vector bundle determined by a single cocycle $\bar{g}_{1,2}: (V_1\cap V_{2}) \times D \rightarrow \GL_n(\CX)$. First note that both $V_1\times D$ and $V_2\times D$ are contractible topological spaces so any  vector bundle over these spaces is topologically trivial.  Furthermore, both of these sets are Stein manifolds (\cite{GunRos}, p.~209) so any topologically trivial analytic vector bundle is analytically trivial (\cite{grauert}, Satz 2, \cite{Oka_theory}, Corollary 3.2). Therefore, $V_1\times D$ and $V_2\times D$ form a covering such that the vector bundle is  analytically trivial on each set.  This implies that the vector bundle is indeed determined by a cocycle $\bar{g}_{1,2}: (V_1\cap V_2) \times D \rightarrow \GL_n(\CX)$.\\[0.1in] 
Under the isomorphism described above, the connections $\nabla(t)$ have a new coordinate description. If we consider the covering  $\{U_i\}_{1\leq i \leq N}$ of $\PX^1(\CX)$, the isomorphism yields equivalent cocycles $\bar{g}_{i,j}(t) = \Gamma_i(t)^{-1} g_{i,j}(t)\Gamma_j(t)$ where  $\Gamma_i(t)$ and~ $\Gamma_j(t)$ are analytically invertible matrices on their respective  coordinate patches.   Since we are dealing with an isomorphism of vector bundles over $\PX^1(\CX)\times D$, one sees that the $\Gamma_i(t)$ are holomorphic in $t$ as well. The local forms of the connection in this new coordinate description become
\begin{eqnarray}
dy & = &\bar{\omega}_i(t)y,
\end{eqnarray}
where
\begin{eqnarray*} 
\bar{\omega}_i(t) & = & d(\Gamma_i(t))\Gamma_i(t)^{-1} + \Gamma_i(t) \omega_i(t) \Gamma_i(t)^{-1}.
\end{eqnarray*}
Note that the $\bar{\omega_i}(t)$ also has, at worst, poles of order $1$ at the $a_j$.\\[0.1in]
%
We now turn to  the solution of the  weak parameterized Riemann-Hilbert Problem.  In the solution of the Riemann-Hilbert Problem as presented in Chapter 3 of \cite{anosov_bolibruch}, the authors appeal to the Birkhoff-Grothendieck Theorem. To prove the existence of a system (\ref{eqnRH1}) that solves the parameterized weak Riemann-Hilbert Problem, we shall use a parameterized version of the classical Birkhoff-Grothendieck Theorem. We can assume that $s \geq 2$ (otherwise the Riemann-Hilbert Problem is trivial). Consider the vector bundle $\calF$ on $\PX^1(\CX)\times D$ with respect to the covering $\{V_1\times D, V_2\times D\}$, as defined above, where $p_1 = a_1$ and $p_2 = a_2$. For simplicity of notation we will assume that $a_1 = 0 $ and $a_2 = \infty$.  Let $\bar{g}_{1,2}:(V_1\cap V_2) \times D \rightarrow \GL_n(\CX)$ be the associated cocycle. The parameterized  Birkhoff-Grothendieck Theorem  ({\it cf.} \cite{ma80a}, Proposition 4.1;   \cite{bolibruch02}, Theorem 2;  \cite{bolibruch04}, Theorem A.1) states:

\begin{quotation}\noindent {\it There exists  an open polydisk $D' \subset D$ as well as maps \linebreak
 $\Phi_1:V_1\times D'\rightarrow \GL_n(\CX)$  with $\Phi_1$,  $\Phi_1^{-1}$ analytic on $V_1\times D'$ , and $\Phi_2:V_2\times D'\rightarrow \GL_n(\CX)$, with $\Phi_2$,  $\Phi_2^{-1}$ analytic on $V_2\times D'$ such that 
\[\bar{g}_{1,2} = \Phi_1^{-1}x^{\Lambda}\Phi_2\]
where $\Lambda$ = $\diag(\lambda_1, \ldots ,\lambda_n)$ for integers $\lambda_1\geq\ldots\geq \lambda_n$.}\footnote{One may also select $D'$ to be $D\backslash \Sigma$ where $\Sigma$ is an analytic subset  of codimension 1 and assert that $\Phi_1$ and $\Phi_2$ are meromorphic along  $\Sigma$ but we shall not need this stronger version.}\end{quotation}
We will now consider (an isomorphic copy of) $\calF$ determined by $\Phi_1$ and $\Phi_2$. Note that the cocycle defining this vector bundle is $(x-a_1)^\Lambda:V_1\cap V_2 \times D' \rightarrow \GL_n(\CX)$.   In these new coordinates, for each $t \in D'$, the connection $\nabla(t)$ above $V_1$ corresponds to a linear differential equation
\[\frac{\d Y}{\d x} = A_1(x,t)Y\]
where $A_1$ is analytic in $t$ and analytic in $x$ outside of $\{a_1, a_3, \ldots ,a_s\}$ and having poles of order at most one at these points. Above $V_2$ we have the expression 
\[\frac{\d Y}{\d u} = A_2(u,t)Y\]
where $A_2$ is analytic in $t$ and analytic in $u = 1/x$ outside of $x \in \{a_2, a_3, \ldots ,a_s\}$ and having poles of order at most one at these points. Using the cocycle, we have
\[A_1(x,t) = \frac{\d (x^\Lambda)}{\d x} x^{-\Lambda} + x^{\Lambda}A_2(1/x,t) x^{-\Lambda}.\]
Since the right-hand side of this equation represents a function meromorphic at $x=\infty$, we have that $A_1$ is meromorphic at $\infty$ as well.   Lemma~\ref{palais} implies that the entries of $A_1(x,t)$ are rational in $x$ with coefficients that are functions analytic on $D'$. We therefore have the following solution of the weak parameterized  Riemann- Hilbert Problem.
\begin{theorem}\label{weakRHthm} Let $S=\{a_1, \ldots a_s\}$ be a finite subset of $\PX^1(\CX)$ and $D$ an open  polydisk  in $\CX^r$. Let $\gamma_1, \ldots , \gamma_s$ be generators of $\pi_1(\PX^1(\CX)\backslash S ; a_0)$ for some fixed base-point $a_0\in\PX^1(\CX)\backslash S$,  and  let $M_i:D\rightarrow \GL_n(\CX)$,   $i =1 , \ldots , s$, be analytic maps with $M_1\cdot\ldots\cdot M_s=I_n$. There exists a parameterized linear differential system 
 \begin{eqnarray*}\label{eqnRH1}
 \frac{\d Y}{\d x} &= &A(x,t)Y
 \end{eqnarray*}
with $A\in\gl_n(\calO_{D'}(x))$ for some open polydisk $D' \subset D$,  having only regular singular points, all in $S$,  such that for some parameterized fundamental solution, the parameterized monodromy matrix along each $\gamma_i$ is  $M_i$. Furthermore, given any $a_i \in \{a_1, \ldots, a_s\}$, the entries of $A$ can be chosen to have at worst  simple poles at all $a_j \neq a_i$.
\end{theorem}

As an application of Theorem~\ref{schlesinger} and Theorem~\ref{weakRHthm}, we solve the inverse problem of parameterized Picard-Vessiot Theory in a special case. Let $k$ be a  so-called \linebreak $\Delta_t$-{\em universal field}, that is, a $\Delta_t$-field such that for any $\Delta_t$-field $k_0 \subset k$, $\Delta_t$-finitely generated over $\QX$ and any $\Delta_t$-finitely generated extension $k_1$ of $k_0$, there is a  \linebreak $\Delta_t$-differential $k_0$-isomorphism of $k_1$ into $k$ (see \cite{DAAG}, Chapter III, \S7). Note that $k$ is in particular a differentially closed $\Delta_t$-field. As in Section~\ref{schlsec}, let $k(x)$ denote the $\Delta = \{\d_x, \d_{t_1}, \ldots , \d_{t_r}\}$-field of rational functions in the indeterminate~$x$ with coefficients in $k$, where~$x$ is a $\Delta_t$-constant with $\dx(x)=1$ and $\dx$ commutes with the $\partial_{t_i}$. We shall prove the following
\begin{cor}\label{inverseprop} Let $G$ be a $\Delta_t$-linear differential algebraic group defined over $k$ and assume that $G(k)$ contains a finitely generated subgroup $H$ that is Kolchin-dense in  $G(k)$. Then $G(k)$ is the PPV-group of a PPV-extension of $k(x)$. \end{cor}
\begin{proof} Let $M_1, \ldots, M_s$ generate $H$. We may assume that $M_1\cdot\ldots\cdot M_s = I_n$. Let $k_0$ be a finitely generated $\Delta_t$-field containing the entries of the $M_i$. Theorem~\ref{seidenberg}  (with $\calK = \QX$ and $\calK_1= k_0$) implies that there is an open polydisk $D$ in $\CX^r$ such that we may consider the elements of $k_0$ as meromorphic functions on $D$. By shrinking $D$ if necessary, we may assume that the entries of the $M_i$ are analytic on $D$. Theorem~\ref{weakRHthm} implies that there exists a parameterized linear differential system 
 \begin{eqnarray*}\label{eqnRH1}
 \frac{\d Y}{\d x} &= &A(x,t)Y
 \end{eqnarray*}
with $A\in\gl_n(\calO_{D'}(x))$ for some open polydisk $D' \subset D$,  having only regular singular points, all in $S$, such that for some  fundamental solution, the parameterized monodromy matrix along each $\gamma_i$ is  $M_i$. Let $k_1$ be a differentially, finitely generated  $\Delta_t$-extension of $k_0$ containing the coefficients of powers of $x$ in  the entries of $A$. Since $k$ is $\Delta_t$-universal, we may assume that $k_1 \subset k$. Theorem~\ref{schlesinger} implies that the group generated by the $M_i$ is Kolchin-dense in the PPV-extension corresponding to $\frac{\d Y}{\d x} = A(x,t)Y$ and so this PPV-group must be $G(k)$.\end{proof}
\noindent This proof follows closely the proof in \cite{tretkoff79} where the authors show than any linear algebraic group defined over $\CX$ is a Galois group of a Picard-Vessiot extension. Those authors use a solution of the weak Riemann-Hilbert Problem, Schlesinger's Theorem and the fact that any linear algebraic group contains a Zariski dense finitely generated subgroup.  In contrast, not all linear differential algebraic groups contain dense finitely generated subgroups.  For example,  the Kolchin-closure of any  finitely generated subgroup $H$ of the additive group $\Ga(k)$  is a proper subgroup of $\Ga(k)$. A proof of this in the ordinary case, that is, when  $\Delta_t = \{\d_t\}$ proceeds as follows. Let $H$ be generated by $z_1, \ldots , z_m$ and assume $z_1, \ldots , z_s$ is a basis for the $C$-vector space spanned by $H$. The elements of the group $H$ all satisfy $L(y) = 0$ for $L(y) = wr(y, z_1, \ldots , z_s)$ where $wr( \ldots )$ denotes the wronskian determinant. In \cite{Landesman} and \cite{CaSi} it is furthermore  shown that neither $\Ga(k)$ nor  $\Gm(k)$ is  the PPV-group of any PPV-extension of $k(x)$.  Among linear algebraic groups, this group presents the main obstruction to a linear algebraic group being a PPV-group over $k(z)$ since it is shown in \cite{sin11} that a linear algebraic group $G$ defined over $k$ is a PPV-group of a PPV-extension of $k(x)$ if and only if the identity component of $G$ has no quotient isomorphic to $\Ga(k)$ or $\Gm(k)$. The proof of this latter fact relies on Corollary~\ref{inverseprop}.

\bibliographystyle{amsplain}

\end{document}